\documentclass[conference]{IEEEtran}
\IEEEoverridecommandlockouts
\usepackage{amsmath,amssymb,amsfonts}
\usepackage{algorithmic}
\usepackage{graphicx}
\usepackage{booktabs}
\usepackage{array}
\usepackage{graphicx} 
\usepackage{float}
\usepackage{comment}
\usepackage{xurl}                     
\usepackage[compatibility=false]{caption} 
\usepackage[breaklinks]{hyperref}    
\usepackage{mathtools}               
\allowdisplaybreaks                   
\usepackage{adjustbox}                
\usepackage{graphicx}                 
\usepackage{longtable}                
\usepackage{url} 
\usepackage{caption}  
\usepackage{textcomp}
\usepackage{xcolor}
\usepackage[table]{xcolor} 
\usepackage{booktabs}     
\usepackage{caption}
\usepackage{multirow}      
\usepackage{caption}     
\usepackage{booktabs}
\usepackage{multirow}
\usepackage{tikz}
\usetikzlibrary{arrows.meta,positioning,fit,backgrounds}
\usepackage{booktabs}
\usepackage{longtable}
\usepackage{pdflscape}
\usepackage{array}
\usepackage[T1]{fontenc}
\usepackage{lmodern}
\usepackage{caption}
\usepackage[numbers,sort&compress]{natbib}
\def\BibTeX{{\rm B\kern-.05em{\sc i\kern-.025em b}\kern-.08em
    T\kern-.1667em\lower.7ex\hbox{E}\kern-.125emX}}
\begin{document}

\title{The Train Unit Shunting Problem with Service Scheduling and Workforce Assignment: classical and quantum optimization approaches
}

\author{
\IEEEauthorblockN{Alessia Ciacco\IEEEauthorrefmark{1},
Francesca Guerriero\IEEEauthorrefmark{1},
Eneko Osaba\IEEEauthorrefmark{3}}

\IEEEauthorblockA{\IEEEauthorrefmark{1}Department of Mechanical, Energy and Management Engineering,
University of Calabria, Rende (CS), Italy}

\IEEEauthorblockA{\IEEEauthorrefmark{3}TECNALIA, Basque Research and Technology Alliance (BRTA),
Derio (Bizkaia), Spain}

\IEEEauthorblockA{Emails: alessia.ciacco@unical.it,
francesca.guerriero@unical.it,
eneko.osaba@tecnalia.com}
}

\maketitle

\begin{abstract}
Efficient railway rolling stock management requires the coordinated planning of train-unit assignment, shunting movements, track occupation, service operations, and workforce allocation. These interdependent decisions give rise to a challenging combinatorial optimization problem, particularly as the size of the fleet and the planning horizon increase. This paper introduces the Train Unit Shunting Problem with Service Scheduling and Workforce Assignment, a mixed-integer linear programming formulation that integrates these operational decisions within a unified optimization framework. To improve computational tractability, a preprocessing procedure is developed to reduce the size of the resulting model. The proposed approach is evaluated through computational experiments on instances of increasing size and under different objective-function weight and cost configurations. Classical optimization experiments using Gurobi show that filtering method reduces the number of decision variables and constraints, while its impact on computational time varies across instances. In addition, a constrained quadratic model representation is investigated using D-Wave's hybrid quantum-classical optimization solver.
\end{abstract}

\section{Introduction}

The efficient management of railway rolling stock is essential for passenger railway operations. Activities performed at stations, depots, and shunting yards between train arrival and departure involve interdependent assignment, track allocation, shunting, and service scheduling decisions. Train units must be matched to departing services, positioned on suitable tracks, and undergo any required cleaning, inspection, or maintenance before departure. These operations are closely interconnected, as track assignments and service locations affect accessibility and shunting requirements. 

The Train Unit Shunting Problem (TUSP) was formally introduced by Freling et al.~\cite{freling2005shunting}, who considered the assignment of train units from arriving to departing trains and their parking on available tracks. Subsequent studies progressively extended the problem by incorporating more detailed shunting decisions and railway infrastructure constraints. Kroon et al.~\cite{kroon2008shunting} proposed an integrated approach to passenger train shunting. 

Jacobsen and Pisinger~\cite{jacobsen2011train} jointly considered train parking and maintenance, while Kamenga et al.~\cite{kamenga2019} proposed a formulation combining train parking, maintenance, and capacity management. Van den Broek et al.~\cite{vandenbroeck2022local} addressed the TUSP by jointly optimizing train matching, parking, routing, and service scheduling through a local-search approach. 


The increasing computational complexity of combinatorial optimization and the difficulty of scaling classical methods as problem size and complexity grow have stimulated interest in alternative computing paradigms. Among these, quantum computing has emerged as a promising approach that exploits quantum-mechanical principles to process information in fundamentally different ways from classical computing. Within quantum computing, Quantum Annealing (QA) has recently attracted considerable attention as an approach for solving discrete optimization problems. QA is based on the evolution of a quantum system toward low-energy states of a problem-dependent energy landscape, where low-energy configurations represent high-quality solutions of the corresponding optimization problem \cite{farhi2001quantum}. 
The interest in QA has grown considerably in recent years, with applications covering a wide range of combinatorial optimization problems. These include scheduling~\cite{li2023quantum, ciaccoscuole}, routing~\cite{le2023quantum,osaba2025quantum}, facility location~\cite{ciacco2026facility,malviya2023logistics, ciacco2026school}, and packing~\cite{de2022hybrid,garcia2022comparative}. 

Quantum optimization has also recently been explored in the railway sector. Domino et al.~\cite{domino2022quadratic,domino2023quantum} investigated the application of QA to railway rescheduling and conflict management. 

This paper makes three main contributions. First, the Train Unit Shunting and Service Scheduling Problem with Workforce Assignment (TUSPwSS-WA) is formulated, incorporating 
the following characteristics: track occupation and capacity, shunting operations, service scheduling, workforce assignment, and temporal feasibility within a unified framework.

Second, a preprocessing procedure specifically tailored to the proposed model is developed to eliminate infeasible or unnecessary combinations of decisions and reduce the size of the resulting optimization problem. Third, both classical and quantum-based solution approaches are investigated. The computational evaluation considers instances of increasing size with the time horizon scaled proportionally to the number of trains. The experiments account for different service-cost policies, workforce-cost policies, and cost scenarios, as well as variations in the objective-function weights and shunting costs. To the best of our knowledge, this work provides the first investigation of quantum computing for the newly formulated TUSPwSS-WA.

The remainder of this paper is organized as follows. Section~\ref{sec:problem} describes the considered train unit shunting and service scheduling problem and its operational characteristics. Section~\ref{sec:model} presents the proposed mathematical formulation. 
Section~\ref{subsec:preprocessing_filtering} introduces the preprocessing method developed to reduce the problem size and improve the efficiency of the formulation.
Section~\ref{sec:computational} reports the computational experiments and discusses the impact of preprocessing and the different solution approaches on computational performance and solution quality. Finally, Section~\ref{sec:conclusion} summarizes the main conclusions and discusses directions for future research.

\section{Problem Definition}
\label{sec:problem}

The TUSPwSS-WA addresses the integrated planning of train units from their arrival at a railway facility to their subsequent departure as part of a scheduled 
service. During this interval, each unit must be matched with a compatible departing train, allocated to available tracks, and undergo all mandatory services before departure, with operations coordinated according to workforce availability and qualifications. 

The railway infrastructure comprises a depot tracks of possibly different capacities. Some tracks are dedicated to service activities, while directed arcs represent feasible shunting movements. 
Each train unit remains in its initially assigned depot throughout the planning horizon, which determines the tracks and service facilities available to that unit. The selected depot determines the available tracks and service facilities, while the scheduled departure establishes the deadline for completing all required activities.

Shunting movements relocate units within the depot before or between services, subject to predefined durations and track-capacity limits. Each movement requires a predefined duration, and track occupation is constrained by physical capacity. Mandatory services, such as cleaning, inspection, and maintenance, must be performed on designated service tracks and require predefined processing times. Their feasibility depends on track availability, unit location, and the time remaining before departure.

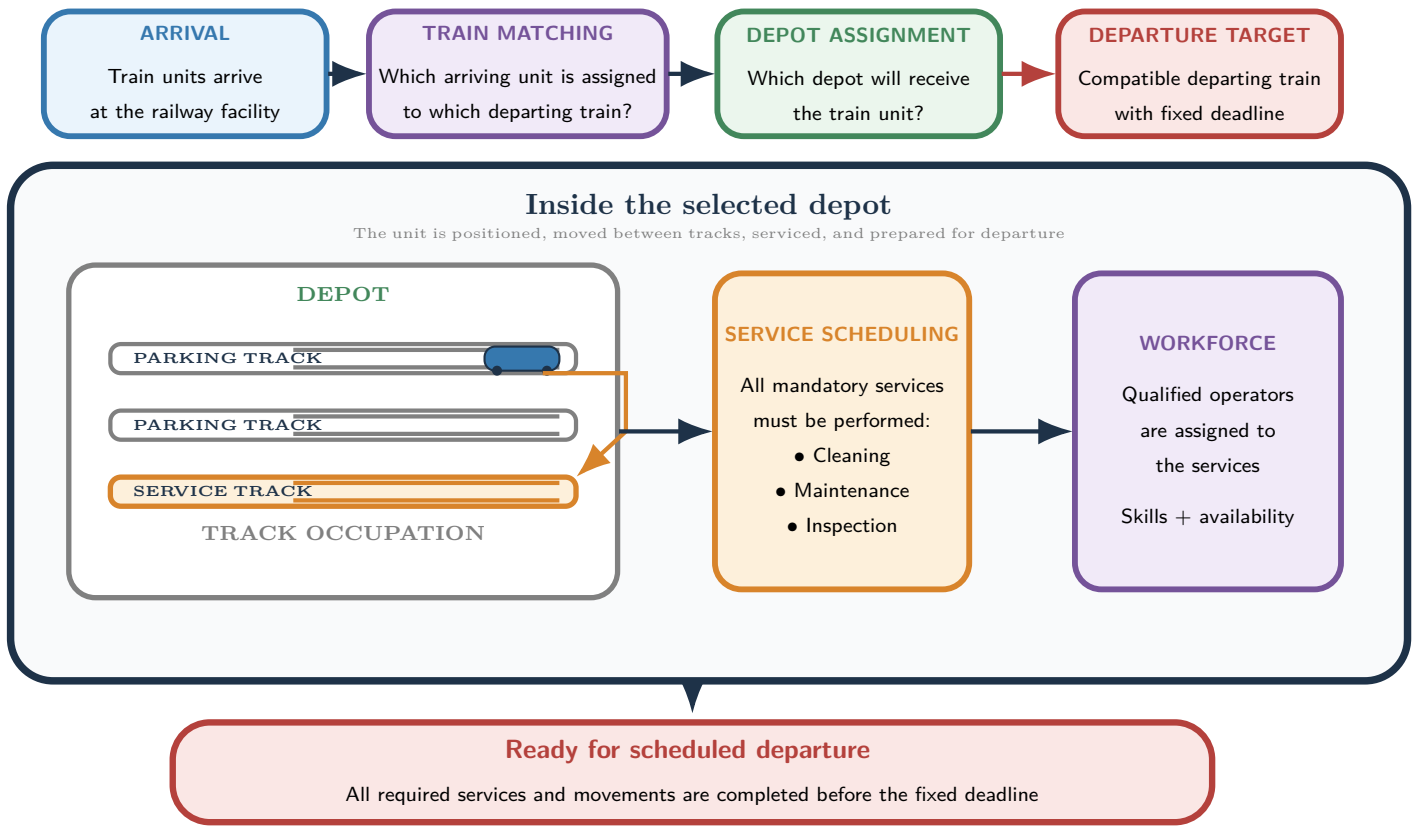
\begin{figure}[ht]
\centering

\begin{tikzpicture}[
    >=Latex,
    font=\sffamily,
    scale=0.55,
    transform shape
]

\definecolor{navy}{RGB}{30,50,72}
\definecolor{blue}{RGB}{54,120,174}
\definecolor{green}{RGB}{66,135,91}
\definecolor{orange}{RGB}{216,132,43}
\definecolor{purple}{RGB}{117,83,153}
\definecolor{red}{RGB}{180,65,60}

\definecolor{lightblue}{RGB}{234,244,252}
\definecolor{lightgreen}{RGB}{235,246,237}
\definecolor{lightorange}{RGB}{253,240,219}
\definecolor{lightpurple}{RGB}{241,234,248}
\definecolor{lightred}{RGB}{250,231,229}
\definecolor{lightgray}{RGB}{249,250,251}


\node[
    draw=blue,
    fill=lightblue,
    rounded corners=5pt,
    minimum width=3.4cm,
    minimum height=1.5cm,
    align=center,
    line width=1.1pt
] (arrival) at (-6.3,0)
{
    {\bfseries\scriptsize\color{blue} ARRIVAL}\\[1mm]
    {\scriptsize Train units arrive}\\
    {\scriptsize at the railway facility}
};

\node[
    draw=purple,
    fill=lightpurple,
    rounded corners=5pt,
    minimum width=3.4cm,
    minimum height=1.5cm,
    align=center,
    line width=1.1pt
] (matching) at (-2.3,0)
{
    {\bfseries\scriptsize\color{purple} TRAIN MATCHING}\\[1mm]
    {\scriptsize Which arriving unit is assigned}\\
    {\scriptsize to which departing train?}
};

\node[
    draw=green,
    fill=lightgreen,
    rounded corners=5pt,
    minimum width=3.4cm,
    minimum height=1.5cm,
    align=center,
    line width=1.1pt
] (depotassign) at (1.8,0)
{
    {\bfseries\scriptsize\color{green} DEPOT ASSIGNMENT}\\[1mm]
    {\scriptsize Which depot will receive}\\
    {\scriptsize the train unit?}
};

\node[
    draw=red,
    fill=lightred,
    rounded corners=5pt,
    minimum width=3.4cm,
    minimum height=1.5cm,
    align=center,
    line width=1.1pt
] (departuretarget) at (5.9,0)
{
    {\bfseries\scriptsize\color{red} DEPARTURE TARGET}\\[1mm]
    {\scriptsize Compatible departing train}\\
    {\scriptsize with fixed deadline}
};

\draw[->, thick, navy] (arrival.east) -- (matching.west);
\draw[->, thick, navy] (matching.east) -- (depotassign.west);
\draw[->, thick, red] (depotassign.east) -- (departuretarget.west);

\node[
    draw=navy,
    fill=lightgray,
    rounded corners=8pt,
    minimum width=16.8cm,
    minimum height=6.2cm,
    line width=1.4pt
] (depotbox) at (0,-4.2) {};

\node[
    font=\bfseries\normalsize,
    text=navy,
    align=center
] at (0,-1.6)
{Inside the selected depot};

\node[
    font=\tiny,
    text=gray,
    align=center
] at (0,-1.93)
{The unit is positioned, moved between tracks, serviced, and prepared for departure};

\node[
    draw=gray,
    fill=white,
    rounded corners=5pt,
    minimum width=6.6cm,
    minimum height=4.0cm,
    line width=1pt
] (tracks) at (-4.4,-4.3) {};

\node[
    font=\bfseries\scriptsize,
    text=green
] at (-4.4,-2.65)
{DEPOT};

\draw[
    draw=gray,
    fill=white,
    rounded corners=2pt,
    line width=0.8pt
] (-7.2,-3.25) rectangle (-1.6,-3.60);

\node[
    anchor=west,
    font=\tiny\bfseries,
    text=navy
] at (-7.05,-3.42)
{PARKING TRACK};

\draw[gray, thick] (-5.0,-3.32) -- (-1.8,-3.32);
\draw[gray, thick] (-5.0,-3.53) -- (-1.8,-3.53);

\draw[
    fill=blue,
    draw=navy,
    rounded corners=2pt
] (-2.7,-3.28) rectangle (-1.8,-3.57);
\fill[navy] (-2.55,-3.57) circle (0.06);
\fill[navy] (-1.95,-3.57) circle (0.06);

\draw[
    draw=gray,
    fill=white,
    rounded corners=2pt,
    line width=0.8pt
] (-7.2,-4.05) rectangle (-1.6,-4.40);

\node[
    anchor=west,
    font=\tiny\bfseries,
    text=navy
] at (-7.05,-4.22)
{PARKING TRACK};

\draw[gray, thick] (-5.0,-4.12) -- (-1.8,-4.12);
\draw[gray, thick] (-5.0,-4.33) -- (-1.8,-4.33);

\draw[
    draw=orange,
    fill=lightorange,
    rounded corners=2pt,
    line width=1pt
] (-7.2,-4.85) rectangle (-1.6,-5.20);

\node[
    anchor=west,
    font=\tiny\bfseries,
    text=navy
] at (-7.05,-5.02)
{SERVICE TRACK};

\draw[orange, thick] (-5.0,-4.92) -- (-1.8,-4.92);
\draw[orange, thick] (-5.0,-5.13) -- (-1.8,-5.13);

\node[
    font=\scriptsize\bfseries,
    text=gray
] at (-4.4,-5.52)
{TRACK OCCUPATION};

\draw[
    ->,
    thick,
    orange
] (-2.0,-3.60)
-- (-1.0,-3.60)
-- (-1.0,-4.30)
-- (-1.6,-4.85);

\node[
    draw=orange,
    fill=lightorange,
    rounded corners=5pt,
    minimum width=0.1cm,
    minimum height=3.8cm,
    align=center,
    line width=1.1pt
] (service) at (1.6,-4.3)
{
    {\bfseries\scriptsize\color{orange} SERVICE SCHEDULING}\\[2mm]
    {\scriptsize All mandatory services}\\
    {\scriptsize must be performed:}\\
    {\scriptsize $\bullet$ Cleaning}\\
    {\scriptsize $\bullet$ Maintenance}\\
    {\scriptsize $\bullet$ Inspection}
};

\node[
    draw=purple,
    fill=lightpurple,
    rounded corners=5pt,
    minimum width=3.2cm,
    minimum height=3.8cm,
    align=center,
    line width=1.1pt
] (workforce) at (6.0,-4.3)
{
    {\bfseries\scriptsize\color{purple} WORKFORCE}\\[2mm]
    {\scriptsize Qualified operators}\\
    {\scriptsize are assigned to}\\
    {\scriptsize the services}\\[2mm]
    {\scriptsize Skills + availability}
};

\draw[->, thick, navy] (tracks.east) -- (service.west);
\draw[->, thick, navy] (service.east) -- (workforce.west);

\draw[
    ->,
    thick,
    navy
] (-0.2,-7.3) -- (-0.2,-7.7);

\node[
    draw=red,
    fill=lightred,
    rounded corners=7pt,
    minimum width=12.5cm,
    minimum height=1.2cm,
    align=center,
    line width=1.2pt
] (ready) at (-0.2,-8.4)
{
    {\bfseries\small\color{red} Ready for scheduled departure }\\[1mm]
    {\scriptsize All required services and movements are completed before the fixed deadline}
};

\end{tikzpicture}

\caption{\small Operational framework of the integrated train unit preparation problem.}
\label{fig:framework}
\end{figure}

As illustrated in Figure~\ref{fig:framework}, the problem jointly considers train matching, track occupation, shunting operations, service scheduling, and workforce allocation.

The workforce is explicitly integrated into the planning problem. Each service requires a predefined number of operators with the appropriate qualifications. Operators cannot perform multiple activities simultaneously, and their total workload is constrained by contractual working-time limits. 

The main decisions are grouped into four interconnected components:
\begin{enumerate}
\item train matching;
\item track occupation and shunting movements;
\item service scheduling; 
\item workforce assignment.
\end{enumerate}

Track assignments determine shunting requirements, movements affect the time available for services, and service schedules determine workforce needs. All activities must respect train arrival times and fixed departure deadlines.

\section{Mathematical Formulation}
\label{sec:model}

The model addresses the integrated operational planning of train units within a railway depot system. In particular, it jointly considers departing-train matching, track occupation, shunting movements, service scheduling, and workforce assignment, assuming that the initial depot of each train unit remains fixed throughout the planning horizon.

\subsection{Sets and Indices}
\begin{itemize}
    \item $U$: set of train units, indexed by $u$
    \item $D$: set of depots, indexed by $d$
    \item $K_d$: set of tracks belonging to depot $d$, indexed by $k$
    \item $K_d^{\mathrm{serv}}\subseteq K_d$: subset of service tracks at depot $d$
    \item $P$: set of departing trains, indexed by $p$
    \item $A_d$: set of directed shunting arcs available at depot $d$, indexed by $a=(k,k')$
    \item $R$: set of service types, indexed by $r$
    \item $R_u\subseteq R$: services required by unit $u$
    \item $W$: set of operators, indexed by $w$
    \item $H$: length of the discrete planning horizon, expressed as the number of time steps
    \item $T=\{0,\ldots,H-1\}$: discrete planning horizon, indexed by $t$
\end{itemize}

\subsection{Parameters}

\begin{itemize}
    \item $a_u$: arrival time of train unit $u$
    \item $l_u$: physical length of train unit $u$
    \item $d(u)$: initial depot of unit $u$, assumed fixed throughout the planning horizon
    \item $k(u)$: initial track of unit $u$
    \item $L_{dk}$: length capacity of track $k$ at depot $d$
    \item $A_p$: scheduled departure time of departing train $p$
    \item $n_p$: number of train units required by departing train $p$
    \item $M_{up} \in \{0, 1\}$: compatibility parameter equal to 1 if train unit $u$ can be assigned to departing train $p$, and 0 otherwise
    \item $\delta_r$: processing duration of service $r$
    \item $\delta_a$: duration of shunting movement $a$
    \item $N_r$: number of operators required to perform service $r$
    \item $Q_{wr} \in \{0, 1\}$: qualification parameter equal to 1 if operator $w$ is qualified to perform service $r$, and 0 otherwise
    \item $H_w^{\max}$: maximum working time allowed for operator $w$
    \item $c^{\mathrm{shunt}}$: unit cost of shunting movements
\item $c^{\mathrm{serv}}_{u,r,k,t}$: cost of performing service $r$ on unit $u$, at depot $d$, on track $k$, starting at time $t$
\item $c^{\mathrm{staff}}_{w,r}$: cost of assigning operator $w$ to service $r$
    \item $w_1, w_2, w_3$: objective-function weights associated with shunting, service and staffing assignment costs, respectively
\end{itemize}

\subsection{Decision Variables}

\begin{itemize}
    \item $m_{up}\in\{0,1\}$: 1 if unit $u$ is matched to departing train $p$
    \item $y_{ukt}\in\{0,1\}$: 1 if unit $u$ occupies track $k$ at time $t$
    \item $q_{uat}\in\{0,1\}$: 1 if train unit $u$ starts shunting movement $a$ at time $t$. A movement started at time $t$
releases its origin track at time $t+1$ and reaches its destination track at time $t+\delta_a$.
    \item $z_{urkt}\in\{0,1\}$: 1 if unit $u$ starts service $r\in R_u$ on service track $k\in K_d^{\mathrm{serv}}$ at time $t$ 
    \item $h_{urwt}\in\{0,1\}$: 1 if operator $w$ is assigned to service $r$ for unit $u$ at time $t$
\end{itemize}

\subsection{Objective Function}

The objective function minimizes a weighted sum of three cost components: shunting movements, service operations, and workforce assignments. It is defined as:
\begin{equation}
\label{eq:objective}
\min \quad w_1 Z^{\text{shunt}} + w_2 Z^{\text{serv}} + w_3 Z^{\text{staff}} \end{equation}

where the individual cost components are:

\begin{flalign}
&Z^{\text{shunt}} = \sum_{u \in U}  \sum_{a \in A_{d(u)}} \sum_{t \in T} c^{\text{shunt}} q_{uat} \label{eq:obj_shunt} \\
&Z^{\text{serv}} =
\sum_{u \in U} \sum_{r \in R_u} 
\sum_{k \in K_{d(u)}^{\text{serv}}} \sum_{t \in T}
c^{\text{serv}}_{u,r,k,t} \;\;z_{urkt}
\label{eq:obj_service} \\
&Z^{\text{staff}} =
\sum_{u \in U} \sum_{r \in R_u} \sum_{w \in W}
\sum_{t \in T}
c^{\text{staff}}_{w,r} h_{urwt}
\label{eq:obj_staff}
\end{flalign}

\noindent \textbf{Shunting Cost Component \eqref{eq:obj_shunt}:} 
    calculates the total cost associated with the movement of train units within the depot network. For every time step $t$, if train unit $u$ initiates shunting arc $a$ at depot $d$ ($q_{uat} = 1$), the unit shunting cost $c^{\text{shunt}}$ is incurred, weighted by $w_1$.
 
\noindent\textbf{Service Operations Cost Component \eqref{eq:obj_service}:}
measures the costs associated with performing the required services ($r \in R_u$) on service tracks ($k \in K_d^{\text{serv}}$). Whenever a service starts at time $t$ ($z_{urkt}=1$), the corresponding cost $c^{\text{serv}}_{u,r,k,t}$ is incurred. This cost depends on the selected service cost policy and its associated parameters, and is weighted by $w_2$.

\noindent\textbf{Workforce Assignment Cost Component \eqref{eq:obj_staff}:}
accounts for the costs of assigning qualified operators to perform services. Whenever operator $w$ is assigned to service $r$ for unit $u$ at time $t$ ($h_{urwt}=1$), the corresponding staff cost $c^{\text{staff}}_{w,r}$ is incurred. This cost depends on the selected staff cost policy and the staff cost multiplier defined by the economic scenario, and is weighted by $w_3$.


\subsection{Constraints}
\vspace{-0.5cm}
\begin{flalign}
&\sum_{u\in U}m_{up}=n_p \quad\forall p\in P, \label{constr:train_demand}\\
&\sum_{p\in P}m_{up}\leq 1 \quad\forall u\in U, \label{constr:unit_matching}\\
&m_{up}\leq M_{up} \quad\forall u\in U,\ p\in P, \label{constr:matching_compatibility}\\
&y_{ukt}=0
\quad\forall u\in U,\ k\in K_{d(u)},\
t\in T:\ t<a_u,
\label{constr:new1}\\
&y_{ukt}=0
\quad\forall u\in U,\
k\in K_{d(u)}\setminus\{k(u)\},\
t\in T:\ t=a_u,
\label{constr:new2}\\
&\sum_{\substack{u\in U\\d(u)=d}}l_u y_{ukt}\leq L_{dk} \quad\forall d\in D,\ k\in K_d,\ t\in T, \label{constr:track_capacity}\\
&q_{u,(k,k'),t}\leq y_{ukt} \quad\forall u\in U,\  (k,k')\in A_{d(u)},\ t\in T, \label{constr:shunt_origin}\\
&y_{uk,t+1}
=
y_{ukt}
-\sum_{k':(k,k')\in A_{d(u)}}q_{u,(k,k'),t}
\nonumber\\
&\quad+
\sum_{\substack{k':(k',k)\in A_{d(u)}\\
t+1-\delta_{(k',k)}\geq a_u}}
q_{u,(k',k),\,t+1-\delta_{(k',k)}}
\nonumber\\
&\quad\forall u\in U,\ k\in K_{d(u)},\
t=a_u,\ldots,H-2.
\label{constr:track_flow}\\
&\sum_{a\in A_{d(u)}} \sum_{\substack{s\in T:\\ s\le t<s+\delta_a}} q_{uas}\le 1 \quad
\forall u\in U,\ t\in T,\label{constr:movimenti_sovrapposti} \\
&z_{urkt}=0
\quad\forall u\in U,\ r\in R_u,\
k\in K_{d(u)}^{\mathrm{serv}},\ t\in T:
t+\delta_r>H,
\label{constr:service_horizon}\\
&q_{uat}=0 \quad \forall u\in U,\ a\in A_{d(u)},\ t\in T:\ t+\delta_a>H-1, \label{constr:orizzonte_temporale} \\
&\sum_{k\in K_{d(u)}^{\mathrm{serv}}}\sum_{t\in T}z_{urkt}=1
\quad\forall u\in U,\ r\in R_u, \label{constr:service_once}\\
&\sum_{r\in R_u}
\sum_{k\in K_{d(u)}^{\mathrm{serv}}}
\sum_{\substack{s\in T:\\s\leq t<s+\delta_r}}
z_{urks}\leq 1
\quad\forall u\in U,\ t\in T,
\label{constr:service_exclusivity}\\
&z_{urkt}\leq y_{uk\tau}
\quad \forall u\in U,\ r\in R_u,\
k\in K_{d(u)}^{\mathrm{serv}},\nonumber\\
&t\in T:\ t+\delta_r\leq H,\ \tau\in T:\ t\leq\tau<t+\delta_r,
\label{constr:service_track}\\
&\sum_{w\in W}h_{urwt}=N_r\sum_{k\in K_{d(u)}^{\mathrm{serv}}}z_{urkt}
\quad\forall u\in U,\ r\in R_u,\ t\in T, \label{constr:staff_requirement_required}\\
&\sum_{u\in U}\sum_{r\in R_u} \sum_{\substack{s\in T:\\s\leq t<s+\delta_r}} h_{urws}\leq 1 \quad\forall w\in W,\ t\in T,
\label{constr:staff_exclusivity}\\
&\sum_{w\in W}h_{urwt}=0
\quad\forall u\in U,\ r\in R\setminus R_u,\ t\in T, \label{constr:staff_requirement_not_required}\\
&\sum_{u\in U}\sum_{r\in R}\sum_{t\in T} \delta_r h_{urwt}\leq H_w^{\max} \quad\forall w\in W, \label{constr:staff_hours}\\
&h_{urwt}\leq Q_{wr}
\quad\forall u\in U,\ r\in R,\ w\in W,\ t\in T, \label{constr:staff_qualification}\\
&t+\delta_r\leq A_p+H(1-m_{up})+H(1-z_{urkt})
\nonumber\\
&\forall u\in U,\ p\in P,\ r\in R_u,\  k\in K_{d(u)}^{\mathrm{serv}},\ t\in T, \label{constr:service_deadline}\\
&y_{u,k(u),a_u}=1 \quad\forall u\in U, \label{constr:initial_track}\\
&m_{up}\in\{0,1\} \quad\forall u\in U,\ p\in P, \label{constr:binary_m}\\
&y_{ukt}\in\{0,1\} \quad\forall u\in U,\  k\in K_{d(u)},\ t\in T, \label{constr:binary_y}\\
&q_{uat}\in\{0,1\} \quad\forall u\in U,\  a\in A_{d(u)},\ t\in T, \label{constr:binary_q}\\
&z_{urkt}\in\{0,1\} \quad\forall u\in U,\ r\in R_u,\  k\in K_{d(u)}^{\mathrm{serv}},\ t\in T,\label{constr:binary_z}\\
&h_{urwt}\in\{0,1\} \quad\forall u\in U,\ r\in R,\ w\in W,\ t\in T. \label{constr:binary_h}
\end{flalign}

Constraints~\eqref{constr:train_demand} ensure that each departing train receives the required number of train units, while constraints~\eqref{constr:unit_matching} ensure that each train unit is assigned to at most one departing train. Constraints~\eqref{constr:matching_compatibility} enforce compatibility between train units and departing trains.
Constraints~\eqref{constr:new1} and~\eqref{constr:new2} define the initial availability of each train unit by preventing track occupation before its arrival and ensuring that, at the arrival time, the unit occupies only its initial track. Constraints~\eqref{constr:track_capacity} limit the total length of train units occupying each track to its physical capacity at every time step. Constraints~\eqref{constr:shunt_origin} allow a shunting movement to start only from a track occupied by the corresponding train unit. Constraints~\eqref{constr:track_flow} describe the evolution of track occupancy over time by accounting for outgoing and incoming shunting movements. Constraints~\eqref{constr:movimenti_sovrapposti} prevent a train unit from performing overlapping shunting movements.
Constraints~\eqref{constr:service_horizon} prevent services from being scheduled if they cannot be completed within the planning horizon, while constraints~\eqref{constr:orizzonte_temporale} exclude shunting movements that cannot reach their destination within the permitted time horizon. Constraints~\eqref{constr:service_once} require each service required by a train unit to be scheduled exactly once. Constraints~\eqref{constr:service_exclusivity} prevent multiple services from being performed simultaneously on the same train unit. Constraints~\eqref{constr:service_track} ensure that the train unit remains on the selected service track throughout the entire service duration.
Constraints~\eqref{constr:staff_requirement_required} enforce the required number of operators for each service required by a train unit at its scheduled start time. Constraints~\eqref{constr:staff_exclusivity} ensure that each operator performs at most one service at any time, taking service durations into account. Constraints~\eqref{constr:staff_requirement_not_required} prevent operators from being assigned to services that are not required by the corresponding train unit. Constraints~\eqref{constr:staff_hours} limit the total working time of each operator to the specified maximum, while constraints~\eqref{constr:staff_qualification} ensure that operators are assigned only to services for which they are qualified.
Constraints~\eqref{constr:service_deadline} ensure that every required service of a train unit matched to a departing train is completed by that train's scheduled departure time. Constraints~\eqref{constr:initial_track} fix each train unit's initial track at its arrival time.
Finally, constraints~\eqref{constr:binary_m}--\eqref{constr:binary_h} define the binary domains of all decision variables.

\section{Preprocessing method}
\label{subsec:preprocessing_filtering}

To reduce the computational burden associated with the full-scale MILP formulation, a preprocessing procedure is introduced to analyze instance data and prune redundant or unfeasible decision variables prior to optimization.

The preprocessing mechanism evaluates temporal horizons and compatibility rules to construct restricted index sets. Since the initial depot of each train unit is fixed throughout the planning horizon, depot-assignment variables are no longer required. Let $\mathcal{M}^{\mathrm{compat}} = \{(u,p) \in U \times P \mid M_{up} = 1\}$ represent the set of compatible unit--train pairs.

Based on these filtered domains and temporal constraints (such as the arrival time $a_u$ of each unit $u$), the index sets of the decision variables are systematically restricted:
\begin{align}
  & m_{up} \in \{0,1\} \quad \forall (u,p) \in \mathcal{M}^{\mathrm{compat}}, \label{eq:filtered_m} \\
  & y_{ukt} \in \{0,1\} \quad \forall u \in U,\; k \in K_{d(u)},\; t \in T:\; t \geq a_u, \label{eq:filtered_y} \\
  & q_{uat} \in \{0,1\} \quad \forall u \in U,\; a \in A_{d(u)},\; t \in T:\;\nonumber\\
& t \geq a_u,\; t+\delta_a \leq H-1,  \label{eq:filtered_q} \\
  & z_{urkt} \in \{0,1\} \quad \forall u \in U,\; r \in R_u,\; k \in K_{d(u)}^{\mathrm{serv}}, \nonumber\\
  & t \in T:\; t+\delta_r \leq H, \label{eq:filtered_z} \\
  & h_{urwt} \in \{0,1\} \quad \forall u \in U,\; r \in R,\; w \in W,\; t \in T. \label{eq:filtered_h}
\end{align}

In the preprocessed optimization model, the objective function \eqref{eq:objective} and all operational constraints are instantiated over the retained decision-variable index sets defined above. The objective function includes only shunting, service, and workforce assignment costs, since depot assignments are fixed inputs rather than decision variables.

The applied filtering rules remove variable index combinations that do not satisfy temporal boundaries or structural compatibility requirements. These reductions preserve the feasible solutions of the original model provided that only combinations that cannot occur in any feasible solution are eliminated.

\section{Computational study}\label{sec:computational}
The goal of the computational experiments is to evaluate the performance TUSPwSS-WA problem obtain with both classical and quantum optimization. 
All computational experiments were conducted in a 64-bit Linux environment (kernel version 6.6.122+) running on an Intel Xeon CPU at 2.20 GHz, with 2 logical CPUs (1 physical core) and 12.67 GB of RAM, using Python 3.13.16. Classical optimization experiments were performed using Gurobi Optimizer version 13.0.3. 

For quantum optimization, we used D-Wave's constrained quadratic model (CQM) hybrid solver, which extends standard quadratic formulations by natively supporting linear constraints and integrality conditions, making it particularly well-suited for mixed-integer optimization problems of this scale.

Unlike purely quantum approaches that require the problem to be mapped directly onto quantum hardware, the hybrid solver combines classical computational methods with quantum processing to search for high-quality feasible solutions. It iteratively explores the solution space, leveraging classical techniques to manage constraints and refine candidate solutions, while exploiting quantum resources as part of the overall optimization process. D-Wave provides a proprietary hybrid solver, \texttt{LeapHybridCQMSampler}, which combines classical heuristics with quantum-guided exploration to achieve a robust balance between global search and local refinement \cite{osaba2024solving, ciacco2026cutting}.

The CQM formulation was implemented using D-Wave’s Ocean SDK (\texttt{dimod} version 0.12.18 and \texttt{dwave-system} version 1.28.0) \cite{DWaveCQMdocs} and executed on D-Wave’s Advantage2\_system1.4 QPU \cite{bertuzzi2024evaluation, osaba2025d, ciacco2025d}.

A time limit of 5 seconds is imposed on all tested instances for classical approaches to ensure a fair comparison under limited computational resources. The time limit for the quantum approach was set to the minimum value required by D-Wave, as determined using its dedicated function.
This section is organized into two main parts. The first describes the instance generation procedure, detailing the methodology adopted to construct the test instances and the characteristics of the resulting datasets. The second presents the computational results, providing a comprehensive evaluation of the proposed approaches in terms of solution quality, computational performance, and the effectiveness of the different optimization strategies.

\subsection{Instance Generation}
\label{subsec:instance_generation}

To the best of our knowledge, no pre-existing railway benchmark dataset specifically suited to the problem considered in this study is available. 
We therefore generate a set of instances covering different problem sizes and operational settings. The instances range from 5 to 20 departing trains, corresponding to 10–40 train units and planning horizons of 15–60 time steps. For each instance size, we consider eight cost configurations obtained by combining two service-cost policies, two workforce-cost policies, and two economic scenarios. The instance generation procedure and the associated parameters are described in detail below. 

All generated instances are available in \cite{repository} to facilitate reproducibility and further research. Due to space limitations, the full description of the generated instances and cost configurations is provided in that accompanying repository. Readers interested in the complete benchmark specification are referred there for additional details.

\begin{table}[htbp]
    \centering
    \tiny
    \caption{Summary of instance characteristics, railway infrastructure, scheduling parameters, service and workforce policies, and economic scenarios used in the computational experiments.}
    \label{tab:experimental_setup}
    \renewcommand{\arraystretch}{1.15}
    \begin{tabular}{p{1.5cm}|| p{6.5cm}}
        \hline
        \textbf{Category} & \textbf{Parameters} \\
        \hline

        \multicolumn{2}{l}{\textit{\textbf{Instance Size and Planning Horizon}}} \\
        \hline
        Instance size
            &  $|P|\in\{5,\ldots,20\}$
               $|U|=2|P|$, yielding 10--40 units \\
      
            & $H=3|P|$, with $T=\{0,\ldots,H-1\}$, yielding 15--60 time steps \\
        \hline

        \multicolumn{2}{l}{\textit{\textbf{Infrastructure and Fleet}}} \\
        \hline
        Railway infrastructure
            & Two depots, each equipped with four tracks of capacity $800\,\mathrm{m}$\\
        Train unit types
            & Two types of train units, with lengths of $40\,\mathrm{m}$ and
              $60\,\mathrm{m}$ \\
        \hline

        \multicolumn{2}{l}{\textit{\textbf{Workforce}}} \\
        \hline
        Number of operators
            & $|W|=\max\{8,2|U|\}=4|P|$, yielding 20--80 operators\\
        Working time limit
            & Each operator's maximum cumulative working time is
              $H_w^{\max}=3H$ \\
        \hline

        \multicolumn{2}{l}{\textit{\textbf{Service Cost Policies}}} \\
        \hline
        Flat Peak
            & $c_{u,r,k,t}=b_r f_{\mathrm{time}}(t)f_{\mathrm{depot}}(d)$,\\
            &  $f_{\mathrm{time}}(t)=1.3$ for $t\geq H/2$
              and $1$ otherwise\\
           &   $f_{\mathrm{depot}}(d)=1.2$ for Depot A and $1$ otherwise \\
        Dynamic Priority
            & $c_{u,r,k,t}=b_r f_{\mathrm{unit}}(u)f_{\mathrm{track}}(k)$,\\
            &  $f_{\mathrm{unit}}(u)=1.4$ for even-indexed
              units and $1$ otherwise\\
          &$f_{\mathrm{track}}(k)=1.15$ for tracks $T_1$ and $T_5$
              and $1$ otherwise\\
        \hline

        \multicolumn{2}{l}{\textit{\textbf{Workforce Cost Policies}}} \\
        \hline
        Standard Seniority
            & $c_{w,r}=m_s b_{w,r}^{\mathrm{staff}}$\\
             & Senior operators: $b_{\mathrm{cleaning}}^{\mathrm{staff}}=8$,
              $b_{\mathrm{maint.}}^{\mathrm{staff}}=25$\\
             & Junior operators: $b_{\mathrm{cleaning}}^{\mathrm{staff}}=5$,
              $b_{\mathrm{maint.}}^{\mathrm{staff}}=15$ \\
        Shift Overtime
            & $c_{w,r}=m_s b_{w,r}^{\mathrm{staff}}$\\
            &  First half of operators:
              $b_{\mathrm{cleaning}}^{\mathrm{staff}}=12$,
              $b_{\mathrm{maint.}}^{\mathrm{staff}}=30$\\
             & Second half of operators:
              $b_{\mathrm{cleaning}}^{\mathrm{staff}}=6$,
              $b_{\mathrm{maint.}}^{\mathrm{staff}}=20$ \\
        \hline

        \multicolumn{2}{l}{\textit{\textbf{Economic Scenarios}}} \\
        \hline
        Balanced
            & $w_1=1.0$, $w_2=2.0$, $w_3=1.5$
              $c_{\mathrm{shunt}}=10$.
               $b_{\mathrm{cleaning}}=15$,
              $b_{\mathrm{maintenance}}=40$
               $m_s=1.0$. \\
        Service and Staff Intensive
            & $w_1=1.0$, $w_2=5.0$, $w_3=4.0$
               $c_{\mathrm{shunt}}=15$
               $b_{\mathrm{cleaning}}=30$,
              $b_{\mathrm{maintenance}}=80$
               $m_s=1.5$. \\
        \hline
    \end{tabular}
\end{table}

\subsection{Numerical Results}
\label{subsec:results}
This section evaluates the computational performance of the proposed MILP formulation and the preprocessing and framework. The analysis is divided into two parts: results obtained using the classical Gurobi solver and results obtained using D-Wave's hybrid quantum-classical solver.

\subsubsection{Classical Optimization Results}
\label{subsec:classical_optimization_results}

This section evaluates the computational performance of the TUSPwSS-WA by comparing the standard formulation with its filtered counterpart. 
Table~\ref{tab:computational_results} reports the results for two selected configurations among all those evaluated, whereas the results obtained for the other configurations are presented in the repository \cite{repository}. 

\begin{table}[ht]
\centering
\caption{\scriptsize Computational comparison between the standard and preprocessed models across all test instances. The columns report: the problem size, in terms of departing trains ($P$), train units ($U$), planning horizon ($H$), shunting unit cost ($c_{\text{sh}}$), the service cost policy (\textbf{SCP}: FP, Flat Peak), the staff cost policy (\textbf{StCP}: SS, Standard Seniority, \textbf{SO}: Shift Overtime) and the economic scenario (\textbf{S}: BAL, Balanced). For each formulation, the table reports the number of decision variables and constraints (\textbf{Var\&Cons}), the objective value (\textbf{Obj}), the optimality gap (\textbf{Opt GAP}(\%)), and the runtime in seconds (\textbf{Time}). \textbf{Std} and \textbf{Prep} denote the standard and preprocessed formulations, respectively.}
\label{tab:computational_results}
\tiny
\setlength{\tabcolsep}{1.8pt}
\begin{tabular}{ccccccc|cccccccc}
\toprule
\multirow{2}{*}{$P$} &
\multirow{2}{*}{$U$} &
\multirow{2}{*}{$H$} &
\multirow{2}{*}{$c_{\text{sh}}$} &
\multirow{2}{*}{SCP} &
\multirow{2}{*}{StCP} &
\multirow{2}{*}{$S$} &
\multicolumn{2}{c}{Var\&Cons} &
\multicolumn{2}{c}{Obj} &
\multicolumn{2}{c}{Opt GAP} &
\multicolumn{2}{c}{Time (s)} \\
\cmidrule(lr){8-9}
\cmidrule(lr){10-11}
\cmidrule(lr){12-13}
\cmidrule(lr){14-15}
& & & & & & & 
Std & Prep &
Std & Prep &
Std & Prep &
Std & Prep \\
\midrule
  5     & 10    & 15    & 10    & FP    & SS    & BAL   & 9619  & 9366  & 766   & 766   & 0     & 0     & 0.16  & 0.09 \\
    6     & 12    & 18    & 10    & FP    & SS    & BAL   & 16004 & 15654 & 713   & 713   & 0     & 0     & 0.10  & 0.11 \\
    7     & 14    & 21    & 10    & FP    & SS    & BAL   & 24645 & 24213 & 918   & 918   & 0     & 0     & 0.26  & 0.27 \\
    8     & 16    & 24    & 10    & FP    & SS    & BAL   & 35940 & 35389 & 1218  & 1218  & 0     & 0     & 0.24  & 0.22 \\
    9     & 18    & 27    & 10    & FP    & SS    & BAL   & 50247 & 49503 & 1565  & 1565  & 0     & 0     & 0.31  & 0.33 \\
    10    & 20    & 30    & 10    & FP    & SS    & BAL   & 68012 & 66956 & 1737  & 1737  & 0     & 0     & 0.42  & 0.87 \\
    11    & 22    & 33    & 10    & FP    & SS    & BAL   & 89545 & 88425 & 1747  & 1747  & 0     & 0     & 0.60  & 0.63 \\
    12    & 24    & 36    & 10    & FP    & SS    & BAL   & 115228 & 113526 & 1832  & 1832  & 0     & 0     & 0.80  & 1.39 \\
    13    & 26    & 39    & 10    & FP    & SS    & BAL   & 145363 & 143530 & 1844  & 1844  & 0     & 0     & 0.97  & 0.98 \\
    14    & 28    & 42    & 10    & FP    & SS    & BAL   & 180380 & 177674 & 1931  & 1931  & 0     & 0     & 1.59  & 1.28 \\
    15    & 30    & 45    & 10    & FP    & SS    & BAL   & 220537 & 217828 & 1911  & 1911  & 0     & 0     & 2.55  & 1.61 \\
    16    & 32    & 48    & 10    & FP    & SS    & BAL   & 266280 & 262905 & 2006  & 2006  & 0     & 0     & 2.93  & 1.96 \\
    17    & 34    & 51    & 10    & FP    & SS    & BAL   & 317927 & 313556 & 2158  & 2158  & 0     & 0     & 2.19  & 2.57 \\
    18    & 36    & 54    & 10    & FP    & SS    & BAL   & 375800 & 370459 & 2405  & 2405  & 0     & 0     & 2.78  & 2.79 \\
    19    & 38    & 57    & 10    & FP    & SS    & BAL   & 440353 & 433825 & 2545  & 2545  & 0     & 0     & 3.37  & 3.41 \\
    20    & 40    & 60    & 10    & FP    & SS    & BAL   & 511828 & 504408 & 2720  & 2720  & 0     & 0     & 3.96  & 3.76 \\

    \midrule
5     & 10    & 15    & 10    & FP    & SO    & BAL   & 9619  & 9366  & 811   & 811   & 0     & 0     & 0.07  & 0.05 \\
    6     & 12    & 18    & 10    & FP    & SO    & BAL   & 16004 & 15654 & 749   & 749   & 0     & 0     & 0.10  & 0.11 \\
    7     & 14    & 21    & 10    & FP    & SO    & BAL   & 24645 & 24213 & 969   & 969   & 0     & 0     & 0.28  & 0.25 \\
    8     & 16    & 24    & 10    & FP    & SO    & BAL   & 35940 & 35389 & 1290  & 1290  & 0     & 0     & 0.23  & 0.20 \\
    9     & 18    & 27    & 10    & FP    & SO    & BAL   & 50247 & 49503 & 1664  & 1664  & 0     & 0     & 0.32  & 0.34 \\
    10    & 20    & 30    & 10    & FP    & SO    & BAL   & 68012 & 66956 & 1845  & 1845  & 0     & 0     & 0.42  & 0.48 \\
    11    & 22    & 33    & 10    & FP    & SO    & BAL   & 89545 & 88425 & 1852  & 1852  & 0     & 0     & 0.60  & 0.63 \\
    12    & 24    & 36    & 10    & FP    & SO    & BAL   & 115228 & 113526 & 1940  & 1940  & 0     & 0     & 0.79  & 0.80 \\
    13    & 26    & 39    & 10    & FP    & SO    & BAL   & 145363 & 143530 & 1949  & 1949  & 0     & 0     & 0.96  & 1.00 \\
    14    & 28    & 42    & 10    & FP    & SO    & BAL   & 180380 & 177674 & 2039  & 2039  & 0     & 0     & 2.04  & 1.28 \\
    15    & 30    & 45    & 10    & FP    & SO    & BAL   & 220537 & 217828 & 2016  & 2016  & 0     & 0     & 1.59  & 1.53 \\
    16    & 32    & 48    & 10    & FP    & SO    & BAL   & 266280 & 262905 & 2114  & 2114  & 0     & 0     & 2.34  & 1.89 \\
    17    & 34    & 51    & 10    & FP    & SO    & BAL   & 317927 & 313556 & 2275  & 2275  & 0     & 0     & 3.30  & 2.25 \\
    18    & 36    & 54    & 10    & FP    & SO    & BAL   & 375800 & 370459 & 2537  & 2537  & 0     & 0     & 2.86  & 2.95 \\
    19    & 38    & 57    & 10    & FP    & SO    & BAL   & 440353 & 433825 & 2957.2 & 2545  & 73    & 0     & 3.99  & 5.04 \\
    20    & 40    & 60    & 10    & FP    & SO    & BAL   & 511828 & 504408 & 3473  & 3414.8 & 85    & 78    & 5.04  & 5.03 \\
   \bottomrule
\tabularnewline
\end{tabular}
\end{table}

The comparison between the Standard Seniority (SS) and Shift Overtime (SO) workforce cost policies, considering the preprocessed formulation, highlights a generally higher objective value under SO. For instances with 5 to 18 departing trains, the increase ranges from 5.05\% to 6.33\%, with an average of approximately 5.71\%, indicating a relatively stable impact of the workforce cost policy across different problem sizes. An exception occurs for (P=19), where both policies yield the same optimal objective value of 2,545. For (P=20), the difference emerges again, with the objective increasing from 2,720 under SS to 2,876 under SO (+5.74\%). Thus, SO generally leads to higher total costs than SS, although this effect is not uniform across all instances.

Regarding computational performance, the standard formulation reaches optimality for all instances under SS, whereas under SO it fails to certify optimality for the two largest instances ((P=19) and (P=20)), exhibiting optimality gaps of 73\% and 85\%, respectively. In contrast, the preprocessed formulation achieves optimality for all tested instances under both workforce policies, with the exception of the largest instance ($P=20$). Nevertheless, even for this instance, preprocessing yields a substantially smaller optimality gap than the standard formulation.

The impact of preprocessing on computational time is less consistent: the average runtime decreases by approximately 4.1\% under SS and 4.4\% under SO, although the preprocessed formulation is not faster for every instance. This behaviour can be partly explained by the relatively short runtimes observed across all instances, for which solver overhead and other computational factors may have a noticeable influence on the overall solution time. Consequently, the reduction in model size does not necessarily translate into a proportional reduction in runtime. Further investigations will be conducted as future work to provide a more comprehensive assessment of the impact of this reduction on solver runtime.

The complete computational results, including all tested instances and the corresponding analysis, are available in the repository \cite{repository}, which readers are encouraged to consult for further details and analysis.

\subsubsection{Quantum Optimization Results}
This section evaluates the computational performance of the TUSPwSS-WA by comparing the standard formulation with its filtered counterpart, using D-Wave's hybrid quantum-classical solver. Table~\ref{tab:quantum_results} summarizes the results obtained for the two configurations presented in Table~\ref{tab:computational_results}. All computational tests were performed over five independent runs for the stochastic nature of the solver and improve the reliability of the performance evaluation by reducing the impact of run-to-run variability.

\begin{table}[ht]
\centering
\caption{\scriptsize Computational comparison between the standard and preprocessed models across all test instances. The columns report: the problem size, in terms of departing trains ($P$), train units ($U$), planning horizon ($H$), shunting unit cost ($c_{\text{sh}}$), the service cost policy (\textbf{SCP}: FP, Flat Peak), the staff cost policy (\textbf{StCP}: SS, Standard Seniority, \textbf{SO}: Shift Overtime) and the economic scenario (\textbf{S}: BAL, Balanced). For each formulation, the table reports the number of decision variables and constraints (\textbf{Var\&Cons}), the objective value (\textbf{Obj}), the optimality gap (\textbf{Opt GAP}(\%)), and the runtime in seconds (\textbf{Time}). \textbf{Std} and \textbf{Prep} denote the standard and preprocessed formulations, respectively.}
\label{tab:quantum_results}
\tiny
\setlength{\tabcolsep}{1.5pt}
\begin{tabular}{ccccccc|cccccccc}
\toprule
\multirow{2}{*}{$P$} &
\multirow{2}{*}{$U$} &
\multirow{2}{*}{$H$} &
\multirow{2}{*}{$c_{\text{sh}}$} &
\multirow{2}{*}{SCP} &
\multirow{2}{*}{StCP} &
\multirow{2}{*}{$S$} &
\multicolumn{2}{c}{Obj avg} &
\multicolumn{2}{c}{Obj dev} &
\multicolumn{2}{c}{Time avg} &
\multicolumn{2}{c}{Time dev} \\
\cmidrule(lr){8-9}
\cmidrule(lr){10-11}
\cmidrule(lr){12-13}
\cmidrule(lr){14-15}
& & & & & & & 
Std & Prep &
Std & Prep &
Std & Prep &
Std & Prep \\
\midrule
      5     & 10    & 15    & 10    & FP    & SS    & BAL   & 766   & 766   & 0.00  & 0.00  & 2.36  & 2.12  & 0.12  & 0.15 \\
    6     & 12    & 15    & 10    & FP    & SS    & BAL   & 713   & 713   & 0.00  & 0.00  & 2.43  & 2.20  & 0.22  & 0.16 \\
      7     & 14    & 15    & 10    & FP    & SS    & BAL   & 918   & 918   & 0.00  & 0.00  & 3.35  & 3.24  & 0.22  & 0.22 \\
    8     & 16    & 24    & 10    & FP    & SS    & BAL   & 1.218 & 1.218 & 0.00  & 0.00  & 4.33  & 3.61  & 0.44  & 0.66 \\
    9     & 18    & 27    & 10    & FP    & SS    & BAL   & 1.565 & 1.565 & 0.00  & 0.00  & 5.00  & 3.97  & 0.87  & 0.43 \\
    10    & 20    & 30    & 10    & FP    & SS    & BAL   & 1.737 & 1.737 & 0.00  & 0.00  & 5.33  & 4.70  & 0.49  & 0.20 \\
    11    & 22    & 33    & 10    & FP    & SS    & BAL   & 1.747 & 1.747 & 0.00  & 0.00  & 4.82  & 4.83  & 0.44  & 0.46 \\
    12    & 24    & 36    & 10    & FP    & SS    & BAL   & 1.832 & 1.832 & 0.00  & 0.00  & 7.84  & 6.74  & 0.45  & 0.55 \\
    13    & 26    & 39    & 10    & FP    & SS    & BAL   & 1.844 & 1.844 & 0.00  & 0.00  & 6.01  & 6.34  & 0.49  & 0.69 \\
    14    & 28    & 42    & 10    & FP    & SS    & BAL   & 1.931 & 1.931 & 0.00  & 0.00  & 8.06  & 7.79  & 0.45  & 0.72 \\
    15    & 30    & 45    & 10    & FP    & SS    & BAL   & 1.911 & 1.911 & 0.00  & 0.00  & 9.51  & 10.90 & 0.75  & 1.91 \\
    16    & 32    & 48    & 10    & FP    & SS    & BAL   & 2.006 & 2.006 & 0.00  & 0.00  & 16.92 & 16.91 & 1.67  & 4.05 \\
        \midrule
5     & 10    & 15    & 10    & FP    & SS    & BAL   & 811   & 811   & 0.00  & 0.00  & 5.13  & 4.80  & 0.47  & 0.53 \\
    6     & 12    & 15    & 10    & FP    & SS    & BAL   & 749   & 749   & 0.00  & 0.00  & 5.66  & 5.62  & 1.87  & 1.74 \\
    7     & 14    & 15    & 10    & FP    & SS    & BAL   & 969   & 969   & 0.00  & 0.00  & 7.25  & 7.32  & 1.56  & 1.07 \\
    8     & 16    & 24    & 0     & FP    & SS    & BAL   & 1290  & 1.290 & 0.00  & 0.00  & 11.52 & 11.97 & 2.19  & 1.78 \\
    9     & 18    & 27    & 0     & FP    & SS    & BAL   & 1664  & 1.664 & 0.00  & 0.00  & 13.00 & 12.18 & 2.90  & 2.74 \\
    10    & 20    & 30    & 0     & FP    & SS    & BAL   & 1845  & 1.845 & 0.00  & 0.00  & 18.14 & 17.54 & 1.70  & 0.95 \\
    11    & 22    & 33    & 0     & FP    & SS    & BAL   & 1852  & 1.852 & 0.00  & 0.00  & 16.77 & 19.11 & 4.34  & 9.27 \\
    12    & 24    & 36    & 0     & FP    & SS    & BAL   & 1940  & 1.940 & 0.00  & 0.00  & 18.56 & 18.77 & 2.08  & 4.53 \\
    13    & 26    & 39    & 0     & FP    & SS    & BAL   & 1949  & 1.949 & 0.00  & 0.00  & 31.90 & 30.88 & 11.70 & 10.32 \\
    14    & 28    & 42    & 0     & FP    & SS    & BAL   & 2039  & 2.039 & 0.00  & 0.00  & 30.02 & 26.62 & 7.70  & 4.01 \\
    15    & 30    & 45    & 0     & FP    & SS    & BAL   & 2016  & 2.016 & 0.00  & 0.00  & 28.19 & 23.88 & 10.51 & 4.74 \\
    16    & 32    & 48    & 10    & FP    & SS    & BAL   & 2226  & 2.114 & 0.00  & 0.00  & 43.35 & 38.89 & 10.72 & 8.04 \\
    17    & 34    & 51    & 10    & FP    & SS    & BAL   & 2397  & 2.397 & 0.00  & 0.00  & 18.80 & 16.88 & 4.09  & 0.66 \\
    18    & 36    & 54    & 10    & FP    & SS    & BAL   & 2685  & 2.685 & 0.00  & 0.00  & 28.82 & 27.83 & 6.57  & 5.22 \\
    19    & 38    & 57    & 10    & FP    & SS    & BAL   & 2844  & 2.724 & 0.00  & 163.77 & 36.71 & 34.03 & 6.03  & 5.22 \\
    20    & 40    & 60    & 0     & FP    & SS    & BAL   & 3268.08 & 2.945 & 200.90 & 94.21 & 30.79 & 34.40 & 9.15  & 10.36 \\ 

   \bottomrule
\tabularnewline
\end{tabular}
\end{table}

The two configurations exhibit substantially different computational behaviours.
For the FP-SS-BAL configuration, the D-Wave hybrid solver consistently achieves solutions matching the Gurobi optimum across all tested instances, for both the standard and preprocessed formulations. Furthermore, preprocessing provides a slight computational advantage, resulting in lower average solution times across the five runs. Since this behaviour remained consistent up to $P=16$, with both formulations reliably attaining optimal solutions, no larger instances were investigated for this configuration.

A substantially different behaviour emerges for the FP-SO-BAL configuration. In this case, the hybrid solver is able to identify optimal solutions, matching the Gurobi optima, for both the standard and preprocessed formulations up to $P=15$. The benefits of the proposed preprocessing approach become particularly evident for the largest instances. For $P=16$, the preprocessed formulation still achieves the Gurobi optimal solution, whereas the standard formulation exhibits an optimality gap of 5\%. For $P=17$ and $P=18$, neither formulation is able to match the Gurobi optimum, with optimality gaps of 5\% and 6\%, respectively.  For $P=19$, the optimality gap decreases from 11\% with the standard formulation to 7\% with preprocessing. An even more significant improvement is observed for $P=20$, where the gap drops from 12\% to just 2\%. These results highlight the contribution of preprocessing in improving solution quality as the problem size increases, particularly when the hybrid solver struggles to identify high-quality solutions.
Furthermore, preprocessing also leads to lower average total solution times compared with the standard formulation.

These findings demonstrate that, while preprocessing provides only marginal improvements in computational time for configurations that are already efficiently solved, its contribution becomes considerably more relevant for harder instances, where it can significantly improve the quality of the solutions found. In the future, 
we will extend the computational experiments to larger values of $P$, in order to investigate whether the benefits of preprocessing are maintained or become even more pronounced as the problem size increases.

\section{Conclusions and Further Work}\label{sec:conclusion}

This paper introduced the TUSPwSS-WA, an integrated optimization framework that jointly addresses train-unit matching, track occupation, shunting operations, service scheduling, and workforce assignment. By incorporating these interdependent operational decisions into a MILP formulation. To improve the computational tractability of the proposed model, a preprocessing procedure was developed to eliminate infeasible and unnecessary combinations of decision variables based on compatibility conditions and temporal constraints. The computational experiments, conducted on instances of increasing size and under different operational and economic configurations, demonstrate that the proposed filtering procedure consistently reduces the number of decision variables and constraints. 
This reduction is important because it results in a more compact formulation, allowing larger problem instances to be solved within reasonable computational times. Furthermore, it improves the solver's ability to find high-quality feasible solutions.

The computational study examined different combinations of service scheduling and workforce assignment policies, together with alternative economic scenarios, to assess the proposed approach under a variety of operating conditions. The results show that the workforce cost policy plays a key role in determining the difficulty of the optimization problem. Under the SS policy, the classical solver consistently identifies optimal solutions, whereas its performance deteriorates on the largest instances when the SO policy is adopted. In this more challenging setting, preprocessing proves particularly effective: it enables the solver to reach optimality for larger instances and substantially reduces the optimality gap in the most difficult case. Overall, these findings indicate that the benefits of preprocessing depend not only on problem size but also on the specific combination of operational policies and cost parameters.

Furthermore, this study provides an investigation into the application of hybrid quantum-classical optimization to the TUSPwSS-WA. The results reveal substantially different behaviours across the two configurations considered. For the FP-SS-BAL configuration, both standard and preprocessed formulations consistently achieve solutions matching the Gurobi optimum, with preprocessing providing a reduction in average solution time. While, for the FP-SO-BAL configuration, the hybrid solver is unable to certify optimality for the larger instances. Nevertheless, preprocessing significantly improves solution quality on the largest instances, reducing the optimality gap with respect to the Gurobi reference solution. This highlights the importance of preprocessing, particularly in the context of quantum computing.

Several future research directions can be identified. First, future work will extend the computational experiments to larger problem instances, with particular emphasis on investigating how the effectiveness of preprocessing evolves as the number of departing trains increases. Second, further investigation of the hybrid quantum-classical approach will be conducted to assess its behaviour on larger instances and under a broader range of operational and economic configurations. In particular, the relationship between preprocessing, problem structure, and solution quality deserves further attention. Finally, additional research could investigate alternative model representations and hybrid optimization strategies to improve computational performance and enable the solution of increasingly realistic railway planning instances.

\color{black}
\section*{Acknowledgments}
Eneko Osaba acknowledges support from the European Union, Q-VERSE project (HORIZON-CL4-2025-03-DIGITAL-EMERGING-02, Grant agreement No: 101293734)
\bibliographystyle{apalike}
{\small
\bibliography{rif}}
\end{document}